%% file: 2016-05-17_formal_duality_on_prime_powers_-_preprint_-_Arxiv.tex
\documentclass[12pt,a4paper]{article}
\usepackage{titlesec}
\titleformat*{\section}{\large\bfseries}
\usepackage{fullpage} 

\input{Header}

\title{\large\textbf{Formal-dual subsets of cyclic groups of prime power order}}
\author{\small Robert Schüler\footnote{Universität Rostock, robert.schueler2@uni-rostock.de}}
\date{\small \today}

\begin{document}
\maketitle

\include*{2016-03-02-Abstract}

\include*{2016-05-17-Sec-Introduction}

\include*{2016-05-17-Sec-Formal-Duality}

\include*{2016-05-17-Sec-Restrictions-on-the-weight-enumerator-function}

\include*{2016-05-17-Sec-Formal-Dual-Sets-odd}

\include*{2016-01-05-Acknowledgement}

\bibliographystyle{plain}
\bibliography{source}

\end{document}

%% file: Header.tex
\usepackage[utf8]{inputenc} 
\usepackage[ngerman,american]{babel} 

\usepackage{mathrsfs} 
\DeclareMathAlphabet{\mathpzc}{OT1}{ptm}{m}{it} 
\usepackage{amsmath, amssymb, amsthm}                                     
\usepackage{mathtools} 

\usepackage{nicefrac} 
\usepackage{faktor} 
\usepackage{paralist} 
\usepackage{enumitem} 
\usepackage{newclude} 

\theoremstyle{plain}
\newtheorem{theorem}{Theorem}[section]
\newtheorem{corollary}[theorem]{Corollary}
\newtheorem{lemma}[theorem]{Lemma}

\newtheorem{example}[theorem]{Example}
\newtheorem{lemmadef}[theorem]{Lemma/Definition}

\theoremstyle{definition}
\newtheorem{definition}[theorem]{Definition}

\newcommand{\NN}{\mathbb{N}} 
\newcommand{\ZZ}{\mathbb{Z}} 
\newcommand{\QQ}{\mathbb{Q}} 
\newcommand{\CC}{\mathbb{C}} 

\newcommand{\scalar}[1]{\left< #1 \right>} 

\DeclareMathOperator{\Gal}{Gal} 
\DeclareMathOperator{\TITO}{TITO} 


%% file: 2016-03-02-Abstract.tex
\abstract{
We study the notion of formal-duality over finite cyclic groups of prime power order as introduced by Cohn, Kumar, Reiher and Schürmann. We will prove that for any cyclic group of odd prime power order, as well as for any cyclic group of order $2^{2l+1}$, there is no primitive pair of formally-dual subsets. This partially proves a conjecture, made by the priorly mentioned authors, that the only cyclic groups with a pair of primitive formally-dual subsets are $\{0\}$ and $\ZZ/4\ZZ$.
}

%% file: 2016-05-17-Sec-Introduction.tex
\section{Introduction}
This work contains new results for the understanding of formal-dual sets. There is a deep interest in the study of the minimization of potential energies over the space of periodic sets.
During several numerical experiments, Cohn, Kumar and Schürmann found in \cite{cohn2009:ground_states_and_formal_duality_relations_in_the_gaussian_core_model} that, in low dimensions, each (at least numerically) optimal configuration has a formally-dual periodic set. This is remarkable since this property is very rare among general periodic-sets. Formal-duality arises from a generalisation of the Poisson-summation formula for lattices and was introduced in \cite[Chapter 6]{cohn2009:ground_states_and_formal_duality_relations_in_the_gaussian_core_model}. Nevertheless, the exact relationship between local optimality of energy minimization problems and formal-dual sets is not well understood so far. Also the characterization of formal-dual periodic sets is incomplete, even in dimension one. The work of Cohn, Kumar, Reiher and Schürmann \cite{Cohn2014:formal_duality_and_generalizations_of_the_poisson_summation_formula} reduced the study of formal-dual periodic sets to the study of formal-dual subsets of finite abelian groups. Furthermore they gave a proof that each abelian group $\ZZ/p^2\ZZ$, where $p$ is an odd prime, does not contain any primitive formal-dual subset (see Definition \ref{def:primitive_subset}). In the case of $\ZZ/4\ZZ$ there is the (up to translation) unique subset $\TITO = \{0,1\}\subset\ZZ/4\ZZ$ which is formally-dual to itself. A conjecture mentioned in \cite[Chapter 4.2]{Cohn2014:formal_duality_and_generalizations_of_the_poisson_summation_formula} claims that $\{0\}\subset \ZZ/\ZZ$ and $\TITO\subset\ZZ/4\ZZ$ are the only primitive formal-dual subsets of cyclic groups.
In this article we will show that there is no primitive formal-dual subset of $\ZZ/p^k\ZZ$ if $p$ is an odd prime. On the way we will also obtain several restrictions for primitive formal-dual subsets of arbitrary finite cyclic groups.\\

 This paper is organized as follows: In Section \ref{sec:Formal_duality} 
we recall the term of formal-duality and clarify all notions that are necessary to understand Theorem \ref{thm:there_is_no_formal_dual_set}. Furthermore we develop an alternative view on formal-duality which is useful for the purposes of this paper. In Section \ref{sec:structure_of_the_weight_enumerator_function} we give some general restrictions on weight enumerators of formal-dual sets (see Definition \ref{def:weight_enumerator}). Using these we are able to prove in Section \ref{sec:formal_dual_sets} the following main theorem of this article:

\begin{theorem}\label{thm:there_is_no_formal_dual_set}
Let $p$ be an odd prime and $k\geq 1$. There is no primitive formal-dual subset of $\ZZ/p^k\ZZ$.
\end{theorem}

%% file: 2016-05-17-Sec-Formal-Duality.tex
\section{Formal duality}\label{sec:Formal_duality}
In this section we clarify all terms to understand the main theorem. In particular, will recall the definition of formal-duality from \cite{Cohn2014:formal_duality_and_generalizations_of_the_poisson_summation_formula} and give an equivalent definition of formal-dual subsets of finite cyclic groups. Furthermore, we recall several results about formal-duality and define the term of a primitive formal-dual set.\\

First, we consider a function which plays a central role in this article.
\begin{definition}\label{def:weight_enumerator}
	Let $T$ be a subset of an (additive) abelian Group $G$. Then the \emph{weight enumerator} of $T$ is defined as
	$$\nu_T\colon G\rightarrow\NN_0, \ \nu_T(y) = \#\{(v,w)\in T\times T \ : \ v-w = y\}.$$
	It counts all differences of elements of $T$ which are equal to $y$.
\end{definition}

Now we recall the definition of formal-dual subsets from \cite[Definition 2.9]{Cohn2014:formal_duality_and_generalizations_of_the_poisson_summation_formula}, and develop an equivalent formulation which relates the weight enumerators directly.
\begin{definition}\label{def:formal_duality_group_and_dual}
	Let $S$ be a subset of the finite abelian group $G$ and $T$ a subset of the dual group $\hat{G}$ (that is the group of all homomorphism $G\rightarrow \CC^\ast$).
	Then $S$ and $T$ are said to be \emph{formally-dual} to each other if
	$$\frac{1}{|T|}\nu_T(y) = \left|\frac{1}{|S|}\sum\limits_{x\in S} \scalar{x,y} \right|^2$$
	for every $y\in \hat{G}$ where $ \scalar{x,y} = y(x)$.
\end{definition}
In \cite[Remark 2.10]{Cohn2014:formal_duality_and_generalizations_of_the_poisson_summation_formula} it was mentioned that this definition is actually symmetric. So one can interchange the roles of $S$ and $T$ in the formula (this is possible via the natural isomorphism between $G$ and $\hat{\hat{G}}$).
For the purposes of this paper it is convenient to use an alternative formulation of formal-duality which relates only the weight enumerators of $S$ and $T$ to each other. 

\begin{lemma}\label{lem:equivalent_formulation_of_formal_duality}
	Let $S$ be a subset of the finite abelian group $G$ and $T$ a subset of the dual group $\hat{G}$. Then $S$ and $T$ are formally-dual to each other if and only if
	$$\frac{|S|^2}{|T|}\nu_T(y) = \sum\limits_{v\in G} \nu_S(v) \scalar{v,y}$$
	for every $y\in \hat{G}$.
\end{lemma}
\begin{proof}
	By a short computation using the facts $\overline{\scalar{v,y}} = \scalar{-v,y}$ and $\scalar{x,y}\cdot\scalar{x^\prime,y} = \scalar{x+x^\prime,y}$, we get
	\begin{align*}
		\left|\frac{1}{|S|}\sum\limits_{x\in S} \scalar{x,y} \right|^2 &= \frac{1}{|S|^2}\left(\sum\limits_{x\in S} \scalar{x,y}\right)\left(\sum\limits_{x\in S} \overline{\scalar{x,y}}\right)\\
		&= \frac{1}{|S|^2}\sum\limits_{x,x^\prime\in S} \scalar{x,y}\overline{\scalar{x^\prime,y}}\\
		&= \frac{1}{|S|^2}\sum\limits_{x,x^\prime\in S} \scalar{x-x^\prime,y}\\
		&= \frac{1}{|S|^2}\sum\limits_{v\in G} \nu_S(v) \scalar{v,y}.
	\end{align*}
	By applying this formula to Definition \ref{def:formal_duality_group_and_dual}, we get the required equivalence.
\end{proof}
In the following we will give a structure analysis of $\nu_S$ on cyclic groups $G=\ZZ/n\ZZ$. On that account we identify $\hat{G}$ with $G$ itself, using $\scalar{x,y} = e^{2\pi i x y / n} = \zeta_{n}^{xy}$ for $\zeta_n \coloneqq e^{2\pi i /n}$. Thereby an equivalent notion of formal-duality arises:

\begin{definition}\label{def:formal_duality_subsets_of_group}
Two sets $S,T\subset \ZZ/n\ZZ$ are said to be \emph{formally-dual} to each other if
$$\frac{|S|^2}{|T|}\nu_T(y) = \sum\limits_{v\in \ZZ/n\ZZ} \nu_S(v) \zeta_{n}^{vy}$$
holds for each $y\in \ZZ/n\ZZ$.
Furthermore, a set $S\subset \ZZ/n\ZZ$ is said to be \emph{formal-dual} if there exists a set $T\subset \ZZ/n\ZZ$ such that $S$ and $T$ are formally-dual to each other.
\end{definition}

Using Lemma \ref{lem:equivalent_formulation_of_formal_duality}, we can reformulate all results about pairs $S\subset G$, $T\subset \hat{G}$ of formally-dual sets as results about pairs of formally-dual sets of $\ZZ/n\ZZ$. One useful example is the following lemma from \cite[end of Proof 2.8]{Cohn2014:formal_duality_and_generalizations_of_the_poisson_summation_formula}:

\begin{lemma}\label{lem:product_of_sizes}
Let $S,T\subset \ZZ/n\ZZ$ be formally-dual to each other. Then
$$n=|S|\cdot |T|.$$
\end{lemma}

Another key definition is the notion of primitive formal-dual subsets, that is:
\begin{definition}\label{def:primitive_subset}
A pair of formally-dual subsets $S,T\subset\ZZ/n\ZZ$ is said to be \emph{primitive} if neither $S$ nor $T$ is contained in some proper coset of $\ZZ/n\ZZ$.
A formal-dual subset $S\subset\ZZ/n\ZZ$ is said to be \emph{primitive} if there is some set $T\subset\ZZ/n\ZZ$ such that $S$ and $T$ form a primitive pair.
\end{definition}
By Lemma \ref{lem:equivalent_formulation_of_formal_duality}, this definition is equivalent to the notion of primitive formal-dual sets in \cite[Section 4.1]{Cohn2014:formal_duality_and_generalizations_of_the_poisson_summation_formula}.
According to \cite[Lemma 4.1, 4.2]{Cohn2014:formal_duality_and_generalizations_of_the_poisson_summation_formula} and the discussion below, any pair of sets that are formally-dual to each other can be constructed from some primitive pair. Therefore it suffices for the characterisation of formal-dual sets to find all primitive formal-dual sets.

We end this section with the currently only known examples of primitive formal-dual sets of finite cyclic groups as stated in \cite[Section 3.1 and discussion before]{Cohn2014:formal_duality_and_generalizations_of_the_poisson_summation_formula}.

\begin{example}\label{ex:primitive_formal_dual_cyclic}
The simplest example is the set $\{0\}\subset\{0\} = \ZZ/\ZZ$ which is formally-dual to itself as is easily checked. This relates to lattices and dual lattices.\\
Another example is the configuration $\TITO = \{0,1\}\subset \ZZ/4\ZZ$ which is formally dual to itself.
In \cite[Section 3.1]{Cohn2014:formal_duality_and_generalizations_of_the_poisson_summation_formula} this was checked by using Definition \ref{def:formal_duality_group_and_dual}, but this is also easily checked by Definition \ref{def:formal_duality_subsets_of_group} after computing the weight enumerator
\begin{center}
\begin{tabular}{r|cccc}
$v$&$0$&$1$&$2$&$3$\\
$\{(x,y) \ : \ x-y=v\}$ & $\{(0,0),(1,1)\}$ & $\{(1,0)\}$ & $\emptyset$ & $\{(0,1)\}$ \\
$\nu_{\TITO}(v)$&$2$&$1$&$0$&$1$
\end{tabular}.
\end{center}
\end{example}

Now we have introduced the vocabulary to understand the main theorem, i.e. Theorem~\ref{thm:there_is_no_formal_dual_set}.

%% file: 2016-05-17-Sec-Restrictions-on-the-weight-enumerator-function.tex
\section{Restrictions of the weight enumerator of formal-dual sets}\label{sec:structure_of_the_weight_enumerator_function}
In this section we will prove some restrictions of the weight enumerators of formal-dual sets. In Theorem \ref{thm:structure_of_nu}, which is the main theorem of this section, we will see that the weight enumerator of a formal-dual subset of some finite cyclic group only takes a small number of different values. Likely, there are similar results about weight enumerators in more general cases.\\

Using the structure of the weight enumerator, we are able to prove several further results. Some of those results are even valid for arbitrary finite cyclic groups.\\
In general, there is no function satisfying all of these conditions as we see in the following section.\\

During this section we will use the notation $\gcd(y,n)$ for $y\in\ZZ/n\ZZ$. Since any representative in some coset of $n\ZZ$ has the same greatest common divisor with $n$, this notation is well defined.
First, we will prove the main theorem of this section:

\begin{theorem}\label{thm:structure_of_nu}
Let $T$ be a formal-dual subset of $\ZZ/n\ZZ$.
Then the weight enumerator satisfies
$$\nu_T(y) = \nu_T(\gcd(y,n))$$
for all $y\in\ZZ/n\ZZ$.
\end{theorem}
\begin{proof}
For all $y\in \ZZ/n\ZZ$ with $\gcd(y,n)=1$ let $\sigma_y\in\Gal (\QQ(\zeta_n):\QQ)$ be defined by $\sigma_y:~\QQ(\zeta_n)\rightarrow \QQ(\zeta_n), \ \zeta_n \mapsto \zeta_n^y$.\\
Since $T$ is a formal-dual set there is 
some set $S\subset\ZZ/n\ZZ$ that is formally-dual to $T$. Then, by definition 
$$\frac{|S|^2}{|T|}\nu_T(y) = \sum\limits_{v\in \ZZ/n\ZZ} \nu_S(v) \zeta_{n}^{vy} = (\ast)$$
for arbitrary $y\in \ZZ/n\ZZ$. Now let $d = \gcd(y,n)$. There is some $y^\prime\in\ZZ/n\ZZ$ relatively prime to $n$ with $y = dy^\prime$. Then
$$\frac{|S|^2}{|T|}\nu_T(y) = \sum\limits_{v\in \ZZ/n\ZZ} \nu_S(v) \zeta_{n}^{vdy^\prime} = \sigma_{y^\prime}\left(\sum\limits_{v\in \ZZ/n\ZZ} \nu_S(v) \zeta_n^{vd}\right)$$
Note that the inner sum is exactly $\frac{|S|^2}{|T|}\nu_T(d)\in\QQ$ since $S$ and $T$ are formally-dual to each other. Since the elements of $\Gal (\QQ(\zeta_n):\QQ)$ fix $\QQ$ point-wise, we have
$$\frac{|S|^2}{|T|}\nu_T(y) = \sigma_{y^\prime} \left(\frac{|S|^2}{|T|}\nu_T(d)\right) = \frac{|S|^2}{|T|}\nu_T(d).$$
This equation is equivalent to the assertion of the theorem since $d=\gcd(y,n)$.
\end{proof}

The strong restriction of \ref{thm:structure_of_nu} can be used to simplify the formulas appearing in Definition \ref{def:formal_duality_subsets_of_group}. In order to do so, we first cite a lemma from \cite[Theorem 271]{Hardy2008:An_introduction_to_the_theory_of_numbers} which helps to compute sums of roots of unity. Then we will slightly generalize it to fit our needs. 

\begin{lemmadef}\label{lem:ramanujan_sum}
For any $n,d\in\NN$ the sum
$$C_n(d) \coloneqq \sum_{\substack{v=1,\dots,n\\ \gcd(v,n)=1}} \zeta_n^{vd}$$
is called \emph{Ramanujan's sum}.
The function $\mu:\NN \rightarrow \{-1,0,1\}$ defined by
$$\mu(m) \coloneqq \begin{cases} (-1)^k \text{ if }m\text{ is square-free, where }k\text{ is the number of prime factors of }m \\ 0 \text{ otherwise}\end{cases}$$
is called \emph{Möbius function}.
Ramanujan's sum can be evaluated as
$$C_n(d) = \sum_{g|\gcd(d,n)} \mu(n/g)g.$$
\end{lemmadef}

\begin{lemmadef}\label{lem:generalised_ramanujan_sum}
For any $n\in\NN$ and $d,e|n$ the following equality holds:
$$C_n(d,e) \coloneqq \sum_{\substack{v\in\ZZ/n\ZZ\\ \gcd(v,n)=e}} \zeta_n^{dv} = \sum_{g|\gcd(d,n/e)} \mu(n/eg)g.$$
Note that for $e=1$ this is exactly Lemma \ref{lem:ramanujan_sum}.
If furthermore $(n/e)|d$, then
$$C_n(d,e) = \varphi(n/e),$$
where $\varphi$ is \emph{Euler's totient function}.
For $d=1$ the stated formula simplifies to
$$C_n(1,e) = \mu(n/e).$$
\end{lemmadef}
\begin{proof}
Using Lemma \ref{lem:ramanujan_sum} the first statement is valid due to the following computation:
\begin{align*}
\sum_{\substack{v\in\ZZ/n\ZZ\\ \gcd(v,n)=e}} \zeta_n^{dv} &= \sum_{\substack{v=1\cdot e,\dots,(n/e)\cdot e \\ \gcd(v,n)=e}} \zeta_n^{dv} = \sum_{\substack{v^\prime=1,\dots,(n/e) \\ \gcd(v^\prime,n/e)=1}} \zeta_{n/e}^{dv^\prime}\\
&= C_{n/e}(d) = \sum_{g|\gcd(d,n/e)} \mu(n/eg)g.
\end{align*}
If $(n/e)|d$ this formula simplifies to
$$C_n(d,e) = \sum_{g|\gcd(d,n/e)} \mu(n/eg)g = \sum_{g|(n/e)} \mu(n/eg)g.$$
Since $\sum_{g|n} \mu(n/g)g = \varphi(n)$, which is a well known identity in number theory (see for example \cite[page 48]{Bundschuh1992:Einfuehrung_in_die_Zahlentheorie}), we have
$$C_n(d,e) = \varphi(n/e).$$
This is the second assertion.
The third assertion easily follows from 
$$C_n(1,e) = \sum_{g|\gcd(1,n/e)} \mu(n/eg)g = \sum_{g=1} \mu(n/eg)g = \mu(n/e).$$
\end{proof}

These lemmata are used for the following corollary which simplifies all sums appearing in Definition \ref{def:formal_duality_subsets_of_group}. They provide a new perspective on formal-duality and are useful to prove non-existence of formal-dual sets under certain conditions.

\begin{corollary}\label{cor:nu_T_simple_form}
Let $S,T\subset\ZZ/n\ZZ$ be formally-dual to each other. Then for each $y\in~\ZZ/n\ZZ$, we have
$$\frac{|S|^2}{|T|} \nu_T(y) = \sum_{e|n} C_n(\gcd(y,n),e)\cdot \nu_S(e),$$
where $C_n(\gcd(y,n),e)$ are natural numbers defined in Lemma \ref{lem:generalised_ramanujan_sum}.
Furthermore, the following special cases turn out to be very useful for our purposes:
\begin{enumerate}
\item $\frac{|S|^2}{|T|} \nu_T(1) = \sum_{e|n} \mu(n/e)\nu_S(e)$ and \label{item:nu_T_1_simple_form}
\item $|S|^2 = \frac{|S|^2}{|T|} \nu_T(0) = \sum_{e|n} \varphi(n/e)\nu_S(e)$. \label{item:nu_T_0_simple_form}
\end{enumerate}
If $n=p^k$ for some prime $p$ this can be simplified further as
\begin{enumerate}[resume]
\item $\frac{|S|^2}{|T|} \nu_T(1) = \nu_S(0) - \nu_S(p^{k-1})$ and \label{item:nu_T_1_simple_form_prime_power}
\item $|S|^2 = \frac{|S|^2}{|T|} \nu_T(0) = \nu_S(0) + \sum_{l=0}^{k-1} (p-1)p^{k-l-1}\nu_S(p^l)$. \label{item:nu_T_0_simple_form_prime_power}
\end{enumerate}
\end{corollary}
\begin{proof}
Let $d = \gcd(y,n)$. Since $T$ is a formal-dual subset of $\ZZ/n\ZZ$ we can apply Theorem \ref{thm:structure_of_nu} to gain
$$\frac{|S|^2}{|T|} \nu_T(y) = \frac{|S|^2}{|T|} \nu_T(d).$$
According to Definition \ref{def:formal_duality_subsets_of_group} with $y=d$ we then have
$$\frac{|S|^2}{|T|} \nu_T(y) = \sum\limits_{v\in \ZZ/n\ZZ} \nu_S(v) \zeta_{p^k}^{dv} = \sum_{e|n} \sum_{\substack{v\in\ZZ/n\ZZ \\ \gcd(v,n) = e}} \nu_S(v) \zeta_{p^k}^{dv}.$$
We again apply Theorem \ref{thm:structure_of_nu} to see that all terms $\nu_S(v)$ appearing in any of the inner sums are in fact the same. This yields
$$\frac{|S|^2}{|T|} \nu_T(y) = \sum_{e|n}  \nu_S(e) \sum_{\substack{v\in\ZZ/n\ZZ \\ \gcd(v,n) = e}} \zeta_{p^k}^{dv}.$$
Now we can apply the definition of $C_n(d,e)$ in Lemma \ref{lem:generalised_ramanujan_sum} to see
$$\frac{|S|^2}{|T|} \nu_T(y) = \sum_{e|n}  \nu_S(e) C_n(d,e).$$
This is exactly the stated assertion.\\

For the special cases \ref{item:nu_T_1_simple_form} and \ref{item:nu_T_1_simple_form_prime_power}, we observe the general formula for $d=1$.
According to Lemma \ref{lem:generalised_ramanujan_sum},
$$\frac{|S|^2}{|T|} \nu_T(1) = \sum_{e|n} C_n(1,e) \nu_S(e) = \sum_{e|n} \mu(n/e)\nu_S(e).$$
If $n=p^k$ we can write all divisors $e$ of $n$ as $p^l$ for $l=0,\dots,k$. Furthermore 
\begin{align*}
\mu(n/p^l) = \mu(p^{k-l}) &= \begin{cases} $1$ &\text{if $p^{k-l}$  is square-free and has an even number of prime divisors}\\ $-1$ &\text{if $p^{k-l}$ is square-free and has an odd number of prime divisors}\\
$0$ &\text{otherwise}\end{cases} \\
&=\begin{cases} $1$ &\text{if $l=k$}\\ $-1$ &\text{if $l=k-1$}\\
$0$ &\text{otherwise}\end{cases}
\end{align*}
Therefore,
$$\frac{|S|^2}{|T|} \nu_T(1) = \sum_{e|n} \mu(n/e)\nu_S(e) = \nu_S(0) - \nu_S(p^{k-1})$$
which is Assertion \ref{item:nu_T_1_simple_form_prime_power}.\\

For the special cases \ref{item:nu_T_0_simple_form} and \ref{item:nu_T_0_simple_form_prime_power} we observe the general formula for $y=0$.
According to Lemma \ref{lem:generalised_ramanujan_sum}, we have
$$\frac{|S|^2}{|T|} \nu_T(0) = \sum_{e|n} C_n(\gcd(0,n),e) \nu_S(e) = \sum_{e|n} C_n(n,e) \nu_S(e) = \sum_{e|n} \varphi(n/e)\nu_S(e)$$
since $(n/e)|n$ for all $e|n$. This proves assertion \ref{item:nu_T_0_simple_form}. To see the validity of assertion \ref{item:nu_T_0_simple_form_prime_power}, we simply write the divisors $e$ as $p^l$ for $l=0,\dots,k$ and use the facts $\varphi(1) = 1$ and $\varphi(p^l) = (p-1)p^{l-1}$ for $l\geq 1$.
\end{proof}

The definition of primitive subsets describes a property of the set $S$ and some set that is formally-dual to $S$. In the following corollary of Theorem \ref{thm:structure_of_nu} we see that in $\ZZ/p^k\ZZ$ all sets that are formally-dual to a primitive set $S$ form a primitive pair of formal-dual sets with $S$. It is not yet confirmed that this is true in general since in general there is no unique set (even up to translations and automorphisms) that is formally-dual to $S$, as can be seen in \cite[Theorem 3.2, Remark 3.3]{Cohn2014:formal_duality_and_generalizations_of_the_poisson_summation_formula}.

\begin{corollary}\label{cor:equivalence_of_primitiv}
Let $S\subset\ZZ/p^k\ZZ$ be a primitive formal-dual set and $T\subset\ZZ/p^k\ZZ$ formally-dual to $S$.
Then $T$ is also primitive. Thus $S$ and $T$ form a primitive pair.
\end{corollary}
\begin{proof}
We recall that by Definition \ref{def:primitive_subset}, $S$ is primitive if it is not contained in a proper coset of $\ZZ/n\ZZ$ and there is some set $T^\prime$ formally-dual to $S$ that also is not contained in a proper coset of $\ZZ/n\ZZ$.

It suffices to show that $T$ is also not contained in a proper coset.\\

First we claim that any formal-dual set $Q\subset\ZZ/p^k\ZZ$ is contained in a proper coset if and only if $\nu_Q(1) = 0$.
Note that $Q$ is contained in a proper coset if and only if all differences of elements of $Q$ are divisible by $p$.
Equivalently, no difference of two elements of $Q$ equals $v\in\ZZ/p^k\ZZ$ with $\gcd(v,p^k)=1$.
Therefore $\nu_Q(v) = 0$ for all $v\in\ZZ/p^k\ZZ$ with $\gcd(v,p^k)=1$.
According to Theorem \ref{thm:structure_of_nu}, this is equivalent to $\nu_Q(1) = 0$.\\

So $T$ is not contained in a proper coset if and only if $\nu_T(1)\neq 0$.
By Definition \ref{def:formal_duality_subsets_of_group}, we have
$$\frac{|S|^2}{|T|}\nu_T(1) = \sum\limits_{v\in \ZZ/p^k\ZZ} \nu_S(v) \zeta_{p^k}^{v} = \frac{|S|^2}{|T^\prime|}\nu_{T^\prime}(1),$$
and according to Lemma \ref{lem:product_of_sizes}, we have $n = |S|\cdot |T| = |S|\cdot |T^\prime|$ and therefore $\nu_T(1) = \nu_{T^\prime}(1)$.
Since $T^\prime$ is not contained in a proper coset, we have $\nu_T(1) = \nu_{T^\prime}(1)\neq 0$.
Altogether, $T$ is not contained in a proper coset.
\end{proof}

At the end of this section we will give some more restrictions of the weight enumerator. It turns out that there are even more structural results for primitive formal-dual subsets of $\ZZ/p^k\ZZ$.

\begin{lemma}\label{lem:weights_of_primitive_sets}
	For an arbitrary subset $S$ of $\ZZ/n\ZZ$ the weight enumerator satisfies $\nu_S(v) \leq |S|$ for all $v\in \ZZ/n\ZZ$.
	If $S$ is a primitive formal-dual subset of $\ZZ/p^k\ZZ$, then $\nu_S(p^{k-1}) < |S|$ and $\nu_S(1) >0$.
	If additionally $|S|=p^l$ for some $l$ with $2l\leq k$ then $\nu_S(1) = 1$.
\end{lemma}
\begin{proof}
	Let $v\in\ZZ/p^k\ZZ$. For every element $x\in S$ there is one unique element $y$ such that $x-y=v$. Therefore $\nu_S(v)$ is also the number of elements $x\in S$ such that $x-v\in S$ which is at most $|S|$.\\
	
	Now we suppose $S$ is primitive as in Definition \ref{def:primitive_subset}. As has been seen in the proof of Corollary \ref{cor:equivalence_of_primitiv} $\nu_S(1) > 0$. 
	Now assume $\nu_S(p^{k-1}) = |S|$. Using Corollary \ref{cor:nu_T_simple_form}.\ref{item:nu_T_1_simple_form_prime_power} we have 
	$$\frac{|S|^2}{|T|} \nu_T(1) = \nu_S(0) - \nu_S(p^{k-1}) = |S| - |S| = 0$$
	where $T\subset \ZZ/p^k\ZZ$ is some set formally-dual to $S$. But then, as discussed in the proof of Corollary \ref{cor:equivalence_of_primitiv}, $T$ would be contained in a proper coset of $\ZZ/p^k\ZZ$. But since $S$ is primitive $T$ has to be primitive too according to Corollary \ref{cor:equivalence_of_primitiv}. This is a contradiction.\\
	
	Finally suppose $S$ is primitive and $|S|=p^l$ with $2l\leq k$. We already discussed $\nu_S(1)~>~0$. Suppose $\nu_S(1) \geq 2 \geq \frac{p}{p-1} \geq \frac{p^{2l-k+1}}{p-1}$.
	Using Corollary \ref{cor:nu_T_simple_form}.\ref{item:nu_T_0_simple_form_prime_power} and\\ $\nu_S(0)~=~|S| > 0$, we obtain the formula
	$$p^{2l} = |S|^2 \geq \nu_S(0) + (p-1)p^{k-1}\nu_S(1) > (p-1)p^{k-1}\frac{p^{2l-k+1}}{p-1} = p^{2l}$$
	which is a contradiction. Therefore, $\nu_S(1) = 1$.
\end{proof}

In the next section we will see that in general there is no function satisfying all of the restrictions given here.

%% file: 2016-05-17-Sec-Formal-Dual-Sets-odd.tex
\section{Formal dual sets of $\ZZ/p^k\ZZ$}\label{sec:formal_dual_sets}

In this section we will use the restrictions developed in Section \ref{sec:structure_of_the_weight_enumerator_function} to show Theorem \ref{thm:there_is_no_formal_dual_set}. First, we show for $\ZZ/p^k\ZZ$ that there is in general no function satisfying all of these restrictions. This instantly yields that, in most cases, there is no primitive formal-dual subset of $\ZZ/p^k\ZZ$. The exceptional case is if we are looking for a formal-dual subset of $\ZZ/p^{2l}\ZZ$ with size $p^l$. For example the functions $\nu:\ZZ/p^{2}\ZZ\rightarrow \NN_0$ defined by $\nu(0) = p$, $\nu(1) = 1$ and $\nu(p) = 0$ satisfy all restrictions of Section \ref{sec:structure_of_the_weight_enumerator_function}. But for odd $p$ there is no subset $S$ of $\ZZ/p^2\ZZ$ with $\nu_S = \nu$. In the subsequent discussion we will see that even in the exceptional case there is no primitive formal-dual set for odd $p$.\\

During this section, we assume without loss of generality that $|S| = p^l$ and $|T| = p^{k-l}$ by using Lemma \ref{lem:product_of_sizes}.\\

The following result is a very strong restriction of possible candidates of primitive formal-dual sets of $\ZZ/p^k\ZZ$, even for $p=2$.
\begin{theorem}\label{thm:restriction_to_special_case}
Let $S\subset\ZZ/p^k\ZZ$ be a primitive formal-dual set. Then $k$ has to be even and $|S| = p^{k/2}$.
\end{theorem}
\begin{proof}
As discussed above, we might assume $|S| = p^l$ and there is a set $T$ formally-dual to $S$ with $|T| = p^{k-l}$.

First assume $2l\geq k$ (and therefore $2(k-l)\leq k$).
Since $S$ and $T$ are formally-dual to each other we can apply Corollary \ref{cor:nu_T_simple_form}.\ref{item:nu_T_1_simple_form_prime_power} to get
\begin{align}
p^{3l-k}\nu_T(1) = \frac{|S|^2}{|T|}\nu_T(1) = \nu_S(0) - \nu_S(p^{k-1}) = p^l - \nu_S(p^{k-1}).\label{form:first_order_restriction_formula}
\end{align}
The set $T$ is a primitive formal-dual set since $S$ is primitive as has been seen in Corollary \ref{cor:equivalence_of_primitiv}.
According to Lemma \ref{lem:weights_of_primitive_sets}, we therefore know $\nu_T(1) = 1$ (note that $|T| = p^{k-l}$ with $2(k-l)\leq k$).
Reordering the terms of (\ref{form:first_order_restriction_formula}) we therefore get
\begin{align}
\nu_S(p^{k-1}) = p^l(1-p^{2l-k}).\label{form:formula_for_p^k-1}
\end{align}
Note that, since $2l\geq k$, the right hand side is a negative number unless $k=2l$ while the left hand side is a non-negative number.
Therefore $k$ is even and $|S| = p^l = p^{k/2}$.\\

Now assume $|S| = p^l$ with $2l\leq k$. Then one can interchange the roles of $S$ and $T$ in the arguments above to get to the same result.
\end{proof}

According to Theorem \ref{thm:restriction_to_special_case}, we only have to study subsets of $\ZZ/p^{2l}\ZZ$ of size $p^l$. In the following we will see, that for odd $p$ none of those subsets is a primitive formal-dual set. Note that all known examples of primitive formal-dual sets of finite cyclic groups as seen in Example \ref{ex:primitive_formal_dual_cyclic} are of the exceptional type.\\

To proof Theorem \ref{thm:special_case_odd} we will need one technical lemma and a very helpful technique which was used first by Gregory Minton in \cite[Lemma 4.4]{Cohn2014:formal_duality_and_generalizations_of_the_poisson_summation_formula}.

\begin{lemma}\label{lemma:equally_coset_distribution}
Let $S\subset\ZZ/n\ZZ$ and $m$ be a divisor of $n$.
If there are exactly $|S|^2/m$ ordered pairs of points of $S$ with a difference divisible by $m$, then each coset of $m\ZZ/n\ZZ$ contains the same number of elements of $S$, namely $\frac{|S|}{m}$. (Note that divisibility by $m|n$ in $\ZZ/n\ZZ$ can be inherited from divisibility by $m$ in $\ZZ$ since this notion is independent from the exact choice of representatives, i.e. if $m|v$ then $m|w$ for all $w \in v + n\ZZ$)
\end{lemma}
\begin{proof}
Two points have a difference divisible by $m$ if and only if they are contained in the same coset of $m\ZZ/n\ZZ$.
Therefore the number of pairs of points in $S$ with a distance divisible by $m$ can also be computed as
$$\sum\limits_{0\leq a \leq m-1} |S\cap (a+m\ZZ/n\ZZ)|^2 = (\ast).$$
It is not difficult to check the following by simple analysis: The only configuration of $m$ variables $\mu_1,\dots,\mu_m\geq 0$ with $\sum_{i=1}^m \mu_i = c$ that minimizes $\sum_{i=1}^m \mu_i^2$ is\\ $\mu_1 = \dots = \mu_m = c/m$.
Especially, since $\sum_{0\leq a \leq m-1} |S\cap (a+m\ZZ/n\ZZ)| = |S|$ this yields
$$ (\ast) \geq \sum\limits_{0\leq a \leq m-1} \left(\frac{|S|}{m}\right)^2 = \frac{|S|^2}{m}.$$
By applying the preconditions this inequality holds with equality and therefore\\
$|S\cap~ (a+~m\ZZ/n\ZZ)| =~ |S|/m$ for all $0\leq a \leq m-1$ which proves the lemma.
\end{proof}

Now we are able to proof the non-existence of formal-dual sets in the exceptional case.

\begin{theorem}\label{thm:special_case_odd}
Let $p$ be an odd prime and $l\geq 1$.
There is no primitive formal-dual subset of $\ZZ/p^{2l}\ZZ$ of size $p^l$.
\end{theorem}
\begin{proof}
Suppose there is such a primitive formal-dual set $S\subset \ZZ/p^{2l}\ZZ$ of size $p^l$. Since $S$ is a primitive formal-dual set we can apply Lemma \ref{lem:weights_of_primitive_sets} to get $\nu_S(1) =1$. Using Theorem \ref{thm:structure_of_nu} we know that the number of ordered pairs with a difference not divisible by $p$ can be computed as
$$\sum\limits_{\substack{v\in\ZZ/p^{2l}\ZZ\\ \gcd(v,p^{2l})=1}} \nu_S(v) = (p-1)p^{2l-1}\nu_S(1) = (p-1)p^{2l-1}.$$
Therefore there are exactly $p^{2l} - (p-1)p^{2l-1} = p^{2l-1} = \frac{|S|^2}{p}$ ordered pairs of points with a distance divisible by $p$.
By Lemma \ref{lemma:equally_coset_distribution} we know that each coset of $p\ZZ/p^{2l}\ZZ$ contains exactly $|S|/p = p^{l-1}$ elements of $S$. Therefore, there are unique values $0\leq a_{i,j} < p^{2l-1}$ such that the elements of $S$ can be denoted as
$$x_{i,j} \equiv i + a_{i,j}\cdot p \mod p^{2l} \text{ for } 0\leq i\leq p-1, 1\leq j \leq p^{l-1}.$$

Furthermore, we can write down the differences of pairs of elements of $S$ lying in $1+p\ZZ/p^{2l}\ZZ$ as
\begin{enumerate}
\item $x_{i,j} - x_{i-1,k} \equiv 1 + (a_{i,j} - a_{i-1,k})p = 1 + A_{i,j,k}\cdot p \mod p^{2l}$ for $1\leq i \leq p-1$, $1\leq j,k \leq p^{l-1}$\label{form:difference_congruent_1_symmetric}
\item $x_{0,j} - x_{p-1,k} \equiv 1 + (a_{0,j} - a_{p-1,k} -1)p = 1 + A_{0,j,k}\cdot p \mod p^{2l}$ for $1\leq j,k \leq p^{l-1}$\label{form:difference_congruent_1_assymetric}
\end{enumerate}
where
$$A_{i,j,k} = a_{i,j} - a_{(i-1 \bmod p), k} - \delta_{i=0} = \begin{cases}a_{i,j} - a_{i-1,k} &\text{ if }i\neq 0\\ a_{0,j} - a_{p-1,k} -1 &\text{ if }i=0\end{cases}.$$
Note that any difference $v\in~1 + p\ZZ/p^{2l}\ZZ$ satisfies $\gcd(v,p^{2l}) = ~1$ and by Theorem \ref{thm:structure_of_nu}, we have $\nu_S(v) = \nu_S(1) = 1$. So for every  $v\in 1 + p\ZZ/p^{2l}\ZZ$ we have exactly one ordered pair of elements of $S$ with difference $v$. This yields that the differences listed in \ref{form:difference_congruent_1_symmetric} and \ref{form:difference_congruent_1_assymetric} are exactly the elements of $1+p\ZZ/p^{2l}\ZZ$. In particular, the coefficients $A_{i,j,k}$ need to be distinct modulo $p^{2l-1}$. Since there are exactly $p\cdot p^{l-1} \cdot p^{l-1} = p^{2l-1}$ coefficients $A_{i,j,k}$, we know that these numbers restricted to $\ZZ/p^{2l-1}\ZZ$ are exactly all elements of $\ZZ/p^{2l-1}\ZZ$. Therefore
\begin{align}
\sum\limits_{\substack{0\leq i\leq p-1\\1\leq j,k \leq p^{l-1}}} A_{i,j,k} \equiv 0 + 1 + 2 + \dots + (p^{2l-1}-1) = \frac{1}{2} p^{2l-1}(p^{2l-1}-1) \mod p^{2l-1}.\label{form:A_ijk_sum_1}
\end{align}
On the other hand, observing the definition of $A_{i,j,k}$, we have
\begin{align}
\sum\limits_{\substack{0\leq i\leq p-1\\1\leq j,k \leq p^{l-1}}} A_{i,j,k} &= \sum\limits_{\substack{0\leq i\leq p-1\\1\leq j,k \leq p^{l-1}}} (a_{i,j} - a_{(i-1 \bmod p), k} - \delta_{i=0})\label{form:A_ijk_sum_2} \\
&= \left(\sum\limits_{\substack{0\leq i\leq p-1\\1\leq j,k \leq p^{l-1}}} a_{i,j}\right) - \left(\sum\limits_{\substack{0\leq i\leq p-1\\1\leq j,k \leq p^{l-1}}} a_{(i-1 \bmod p), k}\right) - p^{2l-2} = -p^{2l-2}. \nonumber
\end{align}
Comparing (\ref{form:A_ijk_sum_1}) and (\ref{form:A_ijk_sum_2}) we obtain
$$-p^{2l-2} \equiv \frac{1}{2} p^{2l-1}(p^{2l-1}-1) \mod p^{2l-1}.$$
For odd $p$ the term $(p^{2l-1}-1)$ is even and therefore
$$\frac{1}{2} p^{2l-1}(p^{2l-1}-1) \equiv 0 \not\equiv -p^{2l-2} \mod p^{2l-1}.$$
This is a contradiction, therefore there is no such set $S$.
Note that for $p=2$ this formula always holds.
\end{proof}

Now it is no difficulty to proof the main theorem of this article. So if we search for more primitive formal-dual sets of cyclic groups of prime power order, we can restrict to subsets of $\ZZ/2^{2l}\ZZ$ of size $2^l$.

\begin{proof}[Proof of Theorem \ref{thm:there_is_no_formal_dual_set}]
Suppose $S$ is a primitive formal-dual set of $\ZZ/p^k\ZZ$, where $p$ is an odd prime and $k\geq 1$. According to Theorem \ref{thm:restriction_to_special_case} $k$ is even and $|S| = p^{k/2}$.
If we write $k=2l$ for $l\geq 1$, we have that $S$ is a primitive subset of $\ZZ/p^{2l}\ZZ$ of size $p^l$ which is impossible due to Theorem \ref{thm:special_case_odd}. This is a contradiction, therefore there is no such set $S$.
\end{proof}

Due to this article we know that primitive formal-dual sets of cyclic groups of prime power order are in some sense very rare. We like to end this paper with a discussion that there also is no primitive formal dual set of $\ZZ/16\ZZ$. To study the still unknown cases $\ZZ/64\ZZ$, $\ZZ/256\ZZ$ and so on, further arguments are needed.

\begin{example}
Suppose $S$ is a primitive formal-dual subset of $\ZZ/16\ZZ$. Then $|S|=4$ as has been seen in Theorem \ref{thm:restriction_to_special_case}.
Following the arguments of the proof of Theorem \ref{thm:special_case_odd} we have $\nu_S(0) = 4$, $\nu_S(1) = 1$ and $\nu_S(8) = 0$. Moreover according to Corollary \ref{cor:nu_T_simple_form}.\ref{item:nu_T_0_simple_form_prime_power} we know
$$16 = |S|^2 = \nu_S(0) + 8\nu_S(1) + 4\nu_S(2) + 2\nu_S(4) + \nu_S(8)$$
which simplifies to the equality
$$4 = 4\nu_S(2) + 2\nu_S(4).$$
Therefore, there are only two possible choices for $\nu_S$ namely
\begin{center}
\begin{tabular}{c|cc}
$d|16$ & $\nu_1(d)$ & $\nu_2(d)$\\
\hline
$1$&$1$&$1$\\
$2$&$1$&$0$\\
$4$&$0$&$2$\\
$8$&$0$&$0$\\
$16\equiv 0$&$4$&$4$\\
\end{tabular}.
\end{center}
In either case it is impossible to construct a set with such weight enumerator. This will be briefly argued. Without loss of generality we can assume $\{0,1\} \subset S$ since there is some pair of points with a difference of $1$ and formal-duality is invariant to translation.
It follows a diagram for $\nu_S = \nu_1$ of $\ZZ/16\ZZ$ where $\bullet$ stands for an element that is assured to be in $S$, $\times$ stands for an element that is surely not in $S$ and $\circ$ is used for unsure elements:
$$\bullet\bullet\times\circ\times\times\circ\circ\times\times\circ\circ\times\times\circ\times.$$
One can explain each $\times$ by the weight enumerator. Elements that provide a difference of $\{1,4,8,12,15\}$ with some $\bullet$ element can not be in $S$ ($\nu_1(4) ~=~ \nu_1(8) ~=~ 0$,  $\nu_1(1) ~=~ \nu_1(15) ~= ~1$).
There has to be a pair of elements with a difference of two since $\nu_1(2) = 1$. Therefore either $3$ or $14$ have to be contained in $S$. Due to the symmetry of the construction we assume without loss of generality $3\in S$ which yields the following diagram:
$$\bullet\bullet\times\bullet\times\times\times\times\times\times\circ\times\times\times\times\times.$$
This is since no element with a difference of $\{1,2,3,4,8,12,13,14,15\}$ to some $\bullet$ element can be in $S$.
Now we see the only choice left is $S = \{0,1,3,10\}$.
But now there are two ordered pairs of points with a difference of $7$ namely $(10,3)$ and $(1,10)$ which is a contradiction to $\nu_1(7) = \nu_1(1) = 1$.\\

Now assume $\nu_S = \nu_2$.
Then we start with the diagram
$$\bullet\bullet\times\times\circ\circ\times\times\times\times\times\times\circ\circ\times\times$$
since no element of $S$ can have a difference of $\{1,2,6,8,10,14,15\}$ to some $\bullet$ element.
Moreover some pair of elements need to have a difference of $3$ since $\nu_2(3) = 1$. Therefore we know $4\in S$ or $13\in S$ but due to the symmetry of the construction we might assume without loss of generality $4\in S$ which yields the following diagram:
$$\bullet\bullet\times\times\bullet\times\times\times\times\times\times\times\times\times\times\times$$
since the distances to some $\bullet$ element can not be $\{1,2,3,6,8,10,13,14,15\}$.
We see that we do not have any possibility left to choose a fourth element, which is again a contradiction.

\end{example}

%% file: 2016-01-05-Acknowledgement.tex
\section*{Acknowledgement}
I like to thank Erik Friese for several inspirations, pointers to literature and verification's of arguments.
Moreover I am grateful to Frieder Ladisch for plenty helpful remarks and inspiring a shorter, more elegant proof.
Furthermore I like to express my sincere gratitude to  Achill Schürmann for suggesting this interesting question as well as many patient reviews of preliminary versions of this paper and for pointing out several mistakes.